\documentclass[11pt]{article}
\usepackage{amsmath, amssymb, amsfonts, amsthm, euscript, txfonts, cite, mathrsfs, relsize, yhmath, graphicx, color}
\usepackage[labelsep=period, font=footnotesize, labelfont=bf]{caption}
\usepackage{dsfont}
\usepackage{comment}

\newtheorem{theorem}{Theorem}[section]
\newtheorem{lem}[theorem]{Lemma}
\newtheorem{proposition}[theorem]{Proposition}

\newtheorem{thm}[theorem]{Theorem}
\newtheorem{definition}[theorem]{{\bf Definition}}
\newtheorem{preremark}[theorem]{{\bf Remark}}
\newenvironment{remark}{\begin{preremark}\rm{\hspace{-0.5 em}{\bf}}}{\end{preremark}}

\DeclareMathAlphabet{\mathpzc}{OT1}{pzc}{m}{it}
\DeclareMathAlphabet{\mathbbmsl}{U}{bbm}{m}{sl}

\renewcommand{\mathsf}{\textsf}

\title{\bf\Large Laplacian integral graphs with a given degree sequence constraint 
\vspace{1cm}}

\author{{\sc\large  A.F.  Novanta{\footnotesize$^{\mathlarger{\bot}}$}, \,  C.S. Oliveira{\footnotesize$^{\mathlarger{\ddag}}$}}, \,    L.S. de Lima{\footnotesize$^{\mathlarger{\dag}}$}\vspace{5mm}
\\{\footnotesize{\sl $^{\mathlarger{\bot}}$Programa de P\'{o}s Gradua\c{c}\~{a}o em Engenharia de Produ\c{c}\~{a}o e Sistemas,}}
\vspace{-1mm}
\\{\footnotesize{\sl Centro Federal de Educa\c{c}\~{a}o Tecnol\'{o}gica Celso Suckow da Fonseca\,(CEFET$/$RJ),}}
\vspace{-1mm}
\\{\footnotesize{\sl CEP 20271-110,  Rio de Janeiro, RJ, Brasil.}}
\vspace{2mm}
\\{\footnotesize{\sl $^{\mathlarger{\bot}}$Departamento de Matem\'{a}tica,}}
\vspace{-1mm}
\\{\footnotesize{\sl Col\'{e}gio Pedro II (CP2),}}
\vspace{-1mm}
\\{\footnotesize{\sl CEP 25025-010,  Duque de Caxias, RJ, Brasil.}}
\vspace{2mm}
\\{\footnotesize{\sl $^{\mathlarger{\ddag}}$Departamento de Matem\'{a}tica,}}
\vspace{-1mm}
\\{\footnotesize{\sl Escola Nacional de Ci\^{e}ncias Estat\'{\i}sticas\,(ENCE-IBGE),}}
\vspace{-1mm}
\\{\footnotesize{\sl CEP 20231-050,  Rio de Janeiro, RJ,   Brasil.}}
\vspace{2mm}
\\{\footnotesize{\sl $^{\mathlarger{\dag}}$Departamento de Administra\c c\~{a}o Geral e Aplicada,}}
\vspace{-1mm}
\\{\footnotesize{\sl Universidade Federal do Paran\'a,(UFPR),}}
\vspace{-1mm}
\\{\footnotesize{\sl CEP 80210-170,  Curitiba, PR, Brasil.}}
\vspace{2mm}
\\\vspace{-6mm}
\\{\footnotesize{
$\text{\itshape{\textsf{anderson\_novanta@yahoo.com.br}}}$ \qquad \qquad
$\text{\itshape{\textsf{carla.oliveira@ibge.gov.br}}}$}\qquad \qquad
$\text{\itshape{\textsf{leonardo.delima@ufpr.br}}}$ \vspace{1.5cm}}}

\date{}

\begin{document}

\maketitle

\begin{abstract}
Let $G$ be a graph on $n$ vertices. The Laplacian matrix of $G$, denoted by $L(G)$, is defined as $L(G)=D(G)-A(G)$, where $A(G)$ is the adjacency matrix of $G$ and $D(G)$ is the diagonal matrix of the vertex degrees of $G$. A graph $G$ is said to be $L$-integral is all eigenvalues of the matrix $L(G)$ are integers. In this paper, we characterize all $L$-integral non-bipartite graphs among all connected graphs with at most two vertices of degree larger than or equal to three.
\end{abstract}

\section{Introduction and Notation}\label{intro} 


Let us consider $G=(V,E)$ as an undirected graph, without loops or multiple edges. Let  $d(G)=(d_1(G), d_2(G), \ldots, d_n(G))$ be the sequence degree of $G$, such that $\Delta(G)=d_1(G) \geq d_2(G)\geq \cdots \geq \delta (G)=d_n(G)$ and $k(G)$ is the vertex connectivity of $G$. We define $\mathcal{G}_1$ as the family of connected graphs with sequence degree in a way that $d_1 \geq 3$ and $1 \leq d_i \leq 2,$ $i = 2, \ldots, n.$ Also, we define  $\mathcal{G}_2$ as the family of connected graphs with sequence degree such that $d_1 \geq d_2 \geq 3$ and $ 1 \leq d_i \leq 2,$ $ i = 3, \ldots, n.$ The \textit{sum} graph, denoted by  $G=G_{1}+G_{2}$, is the graph $G$ such that $V = V_{1} \times V_{2}$ and each pair of vertices $(u_{1}, u_{2})$ and $(v_{1}, v_{2})$ are adjacent in $G$ if and only if $u_{1} = v_{1}$ and $(u_{2}, v_{2})\in E_{2}$ or $u_{2} = v_{2}$ and $(u_{1}, v_{1})\in E_{1}.$  The \textit{Firefly graph}, denoted by $\mathpzc{F}_{r,s,t}$, is the graph with $2r+s+2t+1$ vertices that contain $r$ triangles, $s$ pendant edges and $t$ pendant paths of length 2 sharing a common vertex. We write $A(G)$ for the (0,1)-adjacency matrix of a graph and $D(G)$ for the diagonal matrix of the vertex degrees of $G.$ Also, write $L(G)=D(G)-A(G)$ and $Q(G)=D(G)-A(G)$ for the Laplacian matrix and signless Laplacian matrix of $G$. A graph $G$ is $L$-integral (resp. $Q$-integral) if all of its $L$-eigenvalues (resp. $Q$-eigenvalues) are integers. The spectrum of the Laplacian matrix of $G$ is denoted by $Spec_{L}(G)$ = $\left\{ {\mu_1(G)}^{[n_1]}, {\mu_2(G)}^{[n_2]},  \ldots, {\mu_s(G)}^{[n_s]} \right\},$ where $\mu_i(G)$ is the $i$-th largest Laplacian eigenvalue and $n_i$ is its algebraic multiplicity. The \textit{algebraic connectivity} of $G$ is denoted by $a(G)=\mu_{n-1}(G).$ 

Some infinite families of $L$-integral graphs were characterized in the literature as it can be seen in \cite{kir2,oli1,mer1,mer,kir3}. In particular,  Kirkland in \cite{kir3} determined all  \textit{Laplacian integral} graph such that $\Delta(G) =3$. Motivated by that, we studied all Laplacian integral graphs with at most two vertices of degree greater than or equal to 3. It is worth mentioning that Novanta \emph{et al.} in \cite{nov} found all bipartite $Q$-integral graphs in $\mathcal{G}_1$ and $\mathcal{G}_2.$ Since for bipartite graphs the $L$- and $Q$-eigenvalues coincide all $L$-integral bipartite graphs within those families are solved. In this paper, we determine  all non bipartite $L$-integral graphs in the families $\mathcal{G}_1$ and $\mathcal{G}_2.$ Thus, we state our main result:

\begin{thm}\label{thm:main2}
Let $G$ be a graph on $n\geq9 $ vertices with at most two vertices of degree greater than two. Then $G$ is $L$-integral if and only if $G$ is one of the following: $K_{1,n-1}$, $K_2+K_{1,n-3}$, $K_{2,n-2}$, $\mathpzc{F}_{r,s,0}$, where $s \geq 1$ and $r \geq 1$,  $K_1 \vee (r  K_1 \cup s K_2 \cup K_{1,t})$,  where $t \geq 2$ and $r+s \geq 2$ or $ K_2 \vee (n-2) K_1$.
\end{thm}

The remaining of the paper is organized as follows. In Section 2, we will give some important results that will be needed in the sequel. In Section 3, we present all non bipartite $L$-integral graphs in the family  $\mathcal{G}_1$. In Section 4, we present all non bipartite $L$-integral graphs in the family $\mathcal{G}_2$.

\section{Preliminaries} 
In this section, we present  some results that will be useful to prove the main  results of the paper.

\begin{lem}\cite{fie} \label{conexo}
Let $G$ be a  connected graph. Then $a(G)>0$.
\end{lem}
\begin{lem}\cite{fie} \label{conectividade}
Let $G$ be a non-complete graph. Then $a(G)\leq k(G) \leq \delta (G)$.
\end{lem}

\begin{thm}\cite{kir} \label{criteriodekirkland}
Let $G$ be a non-complete and connected graph on $n$ vertices. Then $k(G)=a(G)$ if only if $G$ can be written as $G=G_{a} \vee G_{b},$ where $G_{a}$ is a disconnected graph on $(n-k(G))$ vertices and $G_{b}$ is a graph on $k(G)$ vertices with $a(G_{b})\ge2k(G)-n$. 
\end{thm}

\begin{lem}\cite{horn} \label{Horn}
Let $A$ be a block diagonal matrix where $A_{ii}$, for $1\leq i \leq k$ are blocks of $A$. Then, $\det \, A=\displaystyle \prod_{i=1}^k \det \, A_{ii}.$  
\end{lem}

\begin{definition}\cite{bro} \label{equitablepartion2}
Given a graph $G$, and a matrix $M=[m_{ij}]$ associated with $G$, a partition $\pi$ of $V(G)$, $V(G)=V_1 \cup \cdots \cup V_k$ is equitable with respect to $G$ and $M$, if for all $i,j \in \{1, 2, \cdots, k\}$ $$\sum_{t \in V_j}m_{st}=d_{ij}  $$
is a constant $d_{ij}$ for any $s\in V_i$.
\end{definition}

\begin{thm} \cite{bro} \label{equitablepartion}
Any eigenvalue of $M_{\pi}$ is also an eigenvalue of $M$.  
\end{thm}

\begin{lem}\cite{heu}\label{interlacearesta}
Let $G$ be a graph on $n$ vertices and $f$ an edge of $G$. If $H \cong G \setminus f$ then
$$\mu_{1}(G) \geq \mu_{1}(H) \geq \mu_{2}(G) \geq \mu_{2}(H) \geq...\geq \mu_{n}(G) \geq \mu_{n}(H).$$
If $H$ is a subgraph of $G$ obtained by removing $r$ edges, then for each $i=1, \ldots, n-r$.
$$\mu_{i}(G) \geq \mu_{i}(H) \geq \mu_{i+r}(G).$$ 
\end{lem}

Let $A$ be a matrix of order $n$ and $1 \leq r \leq n$. The matrix $A_r$ of order $r$ is a  principal submatrix of $A$ obtained by deleting $n-r$ rows and the corresponding columns from $A.$

\begin{proposition}\cite{horn} \label{submatrizprincipal} 
Let $A$ be a Hermitian matrix of order $n,$ let $r$ be an integer with $1 \leq r<n$, and let $A_r$ be a principal submatrix of $A$ of order $r$ with eigenvalues $\lambda_1 \geq \cdots \geq \lambda_n$ and $\theta_1 \geq \cdots \geq \theta_r$ respectively. Then, for each $i=1, \ldots,r$
$$\lambda_i \geq \theta_i \geq \lambda_{i+n-r}.$$
\end{proposition}

\begin{remark} \, Let $B_{n-2}$ be a principal submatrix of $L(G)$. From Proposition \ref{submatrizprincipal}, we have that $\theta_{n-3}(B_{n-2}) \geq \mu_{n-1}(G).$  If $G$ is connected and $B_{n-2}$ has at least two eigenvalues in the interval $(0,1)$, we conclude that $0<\mu_{n-1}(G)<1$.   
\end{remark}

\begin{remark} \cite{nov} \label{caminho}
For $n \geq 7$, there are at least 2 eigenvalues of $L(P_n)$ in the interval $(0, 1)$. 
\end{remark}

\begin{thm} \cite{nov} \label{qintegral}
Let $G$ be a graph on $n$ vertices with at most two vertices of degree greater than or equal to 3. Then $G$ is $Q$-integral if and only if $G$ is one of the following: $K_{1,n-1}$, $K_2+K_{1,n-3}$, $K_{2,n-2}$ or $\Gamma_{1,0,1} \cong P_4[K_2,K_1,K_1,K_2]$.
\end{thm}

The following results characterize cographs from forbidden $P_4$ and show that all cographs are $L$-integral.

\begin{thm} \cite{mer} \label{cograph2}
A graph is \textit{cograph} if and only if it does not have an induced subgraph isomorphic to $P_4$.
\end{thm}

\begin{thm}\cite{mer} \label{cograph1}
If $G$ is  \textit{cograph} then $G$ is a $L$-integral.

\end{thm}

\section{L-integral graphs in $\mathbf{\mathpzc{G}_1}$}
 In \cite{nov}, Novanta \textit{et al.}  characterized all $L$-integral bipartite graphs belonging to $\mathpzc{G}_1$. In this section, we characterize all $L$-integral non-bipartite graphs in $\mathpzc{G}_1$. The graphs that belong to family $\mathpzc{G}_1$ are graphs that contain cycles, paths and pending vertices with one  vertex, say $u$, in common such that $d(u) \geq 3$. Notice that the firefly graphs, $\mathpzc{F}_{r,s,t}$, belong to the family $\mathpzc{G}_1$ such that  $\mathpzc{F}_{0,n-1,0} \cong K_{1,n-1},$ $\mathpzc{F}_{0,0,1} \cong K_{1,2}$ are $L$-integral graphs. Below, we present the main result of this section.   



\begin{thm}\label{teoremaG1}
Let $G \in \mathcal{G}_{1}$ be a graph on $n$ vertices.  Then, $G$ is $L$-integral if only if either $G \cong K_{1,n-1}$ or $G \cong\mathpzc{F}_{r,s,0}$, with $s \geq 1$ and $r \geq 1$.
\end{thm}

\begin{proof}

Let $G \in \mathpzc{G}_{1}$. If $G$ is bipartite, from Theorem \ref{qintegral}, $G$ is $L$-integral if and only if $G \cong K_{1,n-1}$. Now, suppose that $G$ is non-bipartite and $L$-integral. From Lemmas \ref{conexo} and \ref{conectividade}, $0<a(G)\leq k(G)=1$, and consequently $a(G)=1$. From Theorem \ref{criteriodekirkland}, $G \cong G_a \vee G_b$ where $V(G_b)=\{u\}$. As  $G \in \mathpzc{G}_{1}$ and $\forall$ $x \in V(G)$, $d(x)\leq2,$ we have that $G_a \cong r \cdot K_1 \cup s \cdot K_2$ where $r\geq1$ and $s\geq1$. So, $G \cong\mathpzc{F}_{r,s,0}$ which is a cograph. From Theorem \ref{cograph1}, $G$ is $L-$ integral and the result follow.         


\end{proof}

\section{L-integral graphs in $\mathcal{G}_2$}

In \cite{nov}, Novanta \textit{et al.} characterized all $L$-integral bipartite graphs belonging to $\mathpzc{G}_{2}$. Let ${\mathcal{G}^{\prime}_2}$ be the subfamily of non-bipartite graphs belonging to $\mathcal{G}_2$. From Lemmas \ref{conexo} and \ref{conectividade}, we have that $0<a(G)\leq k(G)\leq 2$. So, in order to characterize all $L$-integral graphs in ${\mathcal{G}^{\prime}_2}$ we need to consider the cases: $a(G)=k(G)$ and $a(G)<k(G).$    

    
\noindent{\bf{\textsf{Case 1: $G \in {\mathcal{G}^{\prime}_2}$ and $a(G)=k(G)$. }}} \label{a=k}

\begin{thm} \label{a=k1}
Let $G \in {\mathcal{G}^{\prime}_2}$ with $n \geq 7$ vertices. Then $G$ is $L$-integral if and only if $G \cong K_1 \vee (r\cdot K_1 \cup s\cdot K_2 \cup K_{1,t})$,  where $t \geq 2$ and $r+s \geq 2$ or $G \cong K_2 \vee (n-2)\cdot K_1$.   
\end{thm}

\begin{proof}
 Let $G \in {\mathcal{G}^{\prime}_2}$ with $n \geq 7$ vertices. Suppose that $G$ is $L$-integral. So, $a(G)=k(G)=1$ or $a(G)=k(G)=2$. Firstly, suppose that $a(G)=k(G)=1$. From Theorem \ref{criteriodekirkland}, $G \cong G_a \vee G_b$ where $V(G_b)=\{u\}$ and, consequently, $v \in V(G_a)$. Let $x \in V(G_a)$ that such $x \neq  v$. As $G  \in \mathpzc{G}_2$, $d(x) \leq 2$, we conclude that $G_a \cong r \cdot K_1 \cup s\cdot K_2 \cup K_{1,t}$, where $t \geq 2$ and $r+s \geq 2$. Now, suppose that $a(G)=k(G)=2$. From Theorem \ref{criteriodekirkland}, $G \cong G_a \vee G_b$ such that $G_b$ is a graph on two vertices. So, $u, v \in V(G_b)$ and $G_b \cong K_2$. Let $x \in V(G_a)$. As $G \in {\mathcal{G}^{\prime}_2}$ and $d(x) \leq 2$,  we conclude that $G_a \cong (n-2)\cdot K_1$. Then, $G \cong K_2 \vee (n-2)\cdot K_1$ and the result follows.     
\end{proof}

\noindent{\bf{\textsf{Case 2: $G \in {\mathcal{G}^{\prime}_2}$ and $a(G)<k(G)$. }}} \label{a=k}

In this case, we characterize  all $L-$integral in ${\mathcal{G}^{\prime}_2}$ such that $a(G)<k(G) \leq 2$. If $k(G)=1$, $G$ is not $L$-integral. So, we only need to consider that $k(G)=2$. Then, $G$ has only cycles that contain two vertices $u$ and $v$ of degree larger than or equal to $3$. Consequently, we need to analyze the length of paths with end vertices $u$ and $v$. 

\begin{proposition} \label{maiorque p9}
Let $G \in {\mathcal{G}^{\prime}_2}$ with $n \geq 11$ vertices. If $G$ has a subgraph $P_k$, for $k \geq 9$, with end vertices $u$ and $v$, then $G$ is not $L$-integral.  
\end{proposition}

\begin{proof}
Let $G \in {\mathcal{G}^{\prime}_2}$ with $n \geq 11$ vertices. For $k \geq 9,$ suppose that $G$ contains a path $P_k$ with sequence of vertices $ux_1 \cdots x_{k-2}v$. Let $H$ be the subgraph of $G$ obtained by removing the edges $ux_{1}$ and $x_{k-2}v$. So, $H \cong H_{1} \cup P_{k-2}$, where $H_1$ is a non bipartite graph and $\mu_n(H)=\mu_{n-1}(H)=0$. From Remark \ref{caminho}, $P_{k-2}$ has at least 2 eigenvalues in the interval $(0,1).$ Then, we assume that $0<\mu_{n-2}(H) \leq \mu_{n-3}(H)<1$. From Lemma \ref{interlacearesta}, we conclude that $0<\mu_{n-1}(G) \leq \mu_{n-3}(H)<1$. Therefore, $G$ is not $L$-integral. 
\end{proof}

From Proposition \ref{maiorque p9}, now we need to consider the remaining cases when $G$ has a subgraph $P_k$ for $3 \leq k \leq 8$. First, we consider that $G \in {\mathcal{G}^{\prime}_2}$ is a graph that contains $r$ paths $P^{p}_{n_p}$, for $1 \leq p \leq \Delta(G)=\Delta$  with the sequences of vertices $u{x_1^p} \cdots {x_{n_{p}-2}^p}v$ such that $n_{p} \in \{3, 5, 7 \}$. As $G$ is non bipartite graph, the vertices $u$ and $v$ should be adjacent. By a convenient labeling for the vertices, $L(G)$ can be described in the following way:
$$L(G)=\left[ \begin{array}{ccccc}
     \mathcal{D}_{2 \times2}&\mathcal{T}_{2 \times (n_{1}-2)}&\mathcal{T}_{2\times (n_{2}-2)}&\cdots&\mathcal{T}_{2\times (n_{\Delta}-2)}\\
     {{\mathcal{T}}_{(n_{1}-2)\times 2}}
     &\mathcal{A}_{n_{1}-2}&0_{{(n_{1}-2)}\times {(n_{2}-2)}}&\cdots&{0}_{(n_{1}-2) \times (n_{\Delta}-2)}\\
     {\mathcal{T}}_{(n_{2}-2)\times 2} & {0}_{({n_{2}-2)}\times {(n_{2}-2)}}&\mathcal{A}_{(n_{2}-2)}&\cdots&{0}_{(n_{2} -2) \times (n_{\Delta}-2)}\\
     \vdots&\vdots&\vdots&\ddots&\vdots\\
     {{\mathcal{T}}_{(n_{\Delta}-2)\times 2}}&{0}_{(n_{\Delta}-2) \times (n_{1}-2)}&{0}_{(n_{\Delta}-2) \times (n_{2}-2)}& \cdots&\mathcal{A}_{n_{\Delta}-2}
  
\end{array} \right]$$
where \ $\mathcal{D}=[d_{ij}]_{2 \times 2}$ \ such \ that 
$d_{ij}=\left \{ \begin{array}{ccccc}
   \Delta, & if & i=j&&\\
   -1, & if & i \neq j\\
   
\end{array} \right.$ and \ $\mathcal{T}=[t_{ij}]_{2 \times n_{p-2}}$  such that

$t_{ij} = \left \{ \begin{array}{c}
-1, \ i=j=1 \\
-1, \ i=2 \mbox{ and }  j=n_{p}-2\\
0, \ \mbox{otherwise}. 
\end{array} \right.$

Observe that $A_{n_{p}-2} \in \{ A_1, A_3, A_5\}$, where 
    
${A_{1}}=\left[ \begin{array}{c}
     2     \\
      \end{array} \right]$, ${A_{3}}=\left[ \begin{array}{ccc}
     2 & -1&0   \\
     -1 & 2&-1\\
     0&-1&2
     
\end{array} \right]$ and ${A_{5}}=\left[ \begin{array}{ccccc}
     2 & -1&0&0& 0\\
     -1 & 2&-1&0&0\\
     0&-1&2&-1&0\\
     0&0&-1&2&-1\\
     0&0&0&-1&2\\
 \end{array} \right].$

Notice that since $G$ is non-bipartite, the vertices $u$ and $v$ should be adjacent and we will use this fact to prove Proposition \ref{p5+p7}, \ref{P3+P5} and \ref{P3+P7}.

\begin{proposition}\label{p5+p7}
Let $G \in {\mathcal{G}^{\prime}_2}$ with $n \geq 9$ vertices. If $G$ has at least two subgraphs $P_5$, or at least two subgraphs $P_7$ or one subgraph $P_5$ together with a subgraph $P_7$ with end vertices $u$ and $v$, then $G$ is not $L$-integral. 
\end{proposition}

\begin{proof}
Let $G \in {\mathcal{G}^{\prime}_2}$ with $n \geq 9$ vertices. Suppose that $G$ contains at least two paths $P_5$, or at least two paths $P_7$ or one path $P_5$ together with a path $P_7$  with the sequence of vertices $u{x_1^i} \cdots {x_{n_{p}-2}^i}v$ such that $n_{p} \in \{5,7\}$ and $i \geq 2$. In all cases, let $B_{n-2}$ be the principal submatrix of $L(G)$ obtained by removing both rows and columns that correspond to vertices $u$ and $v.$ It is easy to see that $B_{n-2}$ is a block diagonal matrix and its blocks are the matrices $A_3$ and/or $A_5$ which have eigenvalues in the interval $(0,1)$. From Proposition \ref{submatrizprincipal}, we have $0<\mu_{n-1}(G) \leq \theta_{n-3}(B_{n-2}) < 1$. Then, $G$ is not $L$-integral.

\end{proof}
As $G$ is a non-bipartite graph, the following remaining cases are described as: (i) $G$ has at least one subgraph $P_3$ and one subgraph $P_5$, and (ii) $G$ has at least one subgraph $P_3$ and one subgraph $P_7$. Next, Proposition \eqref{P3+P5} proves case (i), and Proposition \eqref{P3+P7} proves case (ii).

\begin{proposition}\label{P3+P5}
Let $G \in {\mathcal{G}^{\prime}_2}$ with $n \geq 6$ vertices. If $G$ has $s \geq 1$ subgraphs $P_3$ and one subgraph $P_5$ with end vertices $u$ and $v$, then $G$ is not $L$-integral.
\end{proposition}

\begin{proof}
Let $G \in {\mathcal{G}^{\prime}_2}$ with $n \geq 6$ vertices. Recall that $uv \in G$. It is easy to see that $G$ is not a L-integral graph for $s=1$. Suppose that $s \geq 2$. By a convenient labeling of the vertices of $G$, the matrix $L(G)$ can be written in the following way:
$$L(G)=\left[ \begin{array}{cccccccc}
    s+2 & -1 & -1 & \cdots & -1 & -1 & 0 & 0 \\ 
    -1   &s+2& -1 & \cdots & -1 & 0 & 0 & -1 \\
    -1 & -1 & 2 & \cdots & 0 &0 &0 &0 \\ 
    \vdots&\vdots&\vdots & \ddots & \vdots& \vdots & \vdots & \vdots\\
    -1 & -1 & 0 &\cdots & 2 & 0 & 0 & 0 \\
    -1 & 0 & 0 & \cdots & 0 & 2 &  -1 & 0 \\ 
    0 & 0 & 0 & \cdots & 0 & -1 &  2 & -1  \\ 
    0 & -1 & 0 & \cdots & 0 & 0 &  -1 & 2  \\
\end{array} \right].$$

According to Theorem \ref{equitablepartion}, the eigenvalues of the matrix
$$R_{L(G)}=\left[ \begin{array}{cccccc}
   s+2 & -1 & -s &-1 & 0 & 0 \\ 
    -1   &s+2& -s & 0 & 0 & -1 \\
    -1 & -1 & 2 & 0 &0 &0 \\ 
    -1 & 0 & 0 & 2 &  -1 & 0 \\ 
    0 & 0 & 0 & -1 &  2 & -1  \\ 
    0 & -1 & 0 & 0 &  -1 & 2  \\
\end{array} \right]$$ 
are eigenvalues of $L(G)$, whose characteristic polynomial is  $p(\lambda)= \lambda^6+(-2s-12)\lambda^5+(s^2+18s+55)\lambda^4+(-6s^2-56s-120)\lambda^3+(10s^2+70s+125)\lambda^2+(-4s^2-30s-50)\lambda$. As $p(1)=s^2-1>0$, for $s \geq 2$, and $p(2)=-4s<0$, we conclude that there is a root in the interval $(1,2)$, and consequently $G$ is not $L$-integral. 
\end{proof}

\begin{proposition} \label{P3+P7}
Let $G \in {\mathcal{G}^{\prime}_2}$ with $n \geq 8$ vertices. If $G$ has at least $s\geq 1$ subgraphs $P_3$ and one subgraph $P_7$ with end vertices $u$ and $v$,  then $G$ is not $L$-integral.
\end{proposition}

\begin{proof}
Let $G \in {\mathcal{G}^{\prime}_2}$ with $n \geq 8$ vertices. It is easy to see that $G$ is not $L$-integral for $s=1$. Suppose that $s \geq 2$. By a convenient labeling for the vertices, $L(G)$ can be described in the following way:

$$L(G)=\left[ \begin{array}{cccccccccc}
    s+2 & -1 & -1 & \cdots & -1 & -1 & 0 & 0 & 0 & 0 \\ 
    -1   &s+2& -1 & \cdots & -1 & 0 & 0 & 0 & 0 & -1 \\
    -1 & -1 & 2 & \cdots & 0 & 0 & 0 & 0 & 0 & 0 \\ 
    \vdots&\vdots&\vdots & \ddots & \vdots& \vdots & \vdots & \vdots & \vdots & \vdots\\
    -1 & -1 & 0 &\cdots & 2 & 0 & 0 & 0 & 0 & 0 \\
    -1 & 0 & 0 & \cdots & 0 & 2 &  -1 & 0 & 0 & 0 \\ 
    0 & 0 & 0 & \cdots& 0 & -1 &  2 & -1 & 0 & 0 \\ 
    0 & 0 & 0 & \cdots& 0 & 0 &  -1 & 2 & -1 & 0  \\
    0 & 0 & 0 & \cdots& 0 & 0 &  0 & -1 & 2 & -1  \\
    0 & -1 & 0 & \cdots& 0 & 0 &  0 & 0 & -1 & 2  \\
 \end{array} \right].$$
According to Theorem \ref{equitablepartion}, the eigenvalues of the matrix
$$R_{L(G)}=\left[ \begin{array}{cccccccc}
   s+2 & -1 & -s &-1 & 0 & 0& 0 & 0 \\ 
    -1   &s+2& -s & 0 & 0 & 0 & 0 & -1 \\
    -1 & -1 & 2 & 0 & 0 & 0 & 0 & 0 \\ 
    -1 & 0 & 0 & 2 &  -1 & 0& 0 & 0 \\ 
    0 & 0 & 0 & -1 &  2 & -1& 0 & 0  \\ 
    0 & 0 & 0 & 0 &  -1 & 2& -1 & 0  \\
    0 & 0 & 0 & 0 &  0 & -1& 2 & -1  \\
    0 & -1 & 0 & 0 &  0 & 0& -1 & 2  \\

\end{array} \right]$$
are eigenvalues of $L(G)$, whose characteristic polynomial is
$p(\lambda)= \lambda^8+(-2s-16)\lambda^7+(s^2+26s+105)\lambda^6+(-10s^2-132s-364)\lambda^5+(36s^2+330s+714)\lambda^4+(-56s^2-420s-784)\lambda^3+(35s^2+252s+441)\lambda^2 +(-6s^2-56s-98)\lambda$. As $p(3)=-6s+3<0$ and $p(4)=24s^2-32s+8>0$ for $s\geq 2$, we conclude that there is a root in the interval $(3,4)$, and consequently $G$ is not $L$-integral.

\end{proof}

\begin{remark}
If $G \in {\mathcal{G}^{\prime}_2}$ and $G$ has only subgraph $P_3$ with end vertices $u$ and $v$, $a(G)=k(G)$ which was analyzed in Case \ref{a=k1}. 
\end{remark}

Now let us analyze the cases in which $G \in {\mathcal{G}^{\prime}_2}$ is a graph that contains $r$ paths $P^{p}_{n_p}$ with the sequence of vertices $u{x_1^p} \cdots {x_{n_{p}-2}^p}v,$ for $1 \leq p \leq r,$ and $n_{p} \in \{4, 6, 8 \}.$ By a convenient labeling to the vertices of $G$ we have
$$L(G)=\left[ \begin{array}{ccccc}
     \mathcal{D}_{2 \times2}&\mathcal{T}_{2 \times (n_{1}-2)}&\mathcal{T}_{2\times (n_{2}-2)}&\cdots&\mathcal{T}_{2\times (n_{r}-2)}\\
     {{\mathcal{T}}_{(n_{1}-2)\times 2}}
     &\mathcal{A}_{n_{1}-2}&0_{{(n_{1}-2)}\times {(n_{2}-2)}}&\cdots&{0}_{(n_{1}-2) \times (n_{r}-2)}\\
     {\mathcal{T}}_{(n_{2}-2)\times 2} & {0}_{({n_{2}-2)}\times {(n_{2}-2)}}&\mathcal{A}_{(n_{2}-2)}&\cdots&{0}_{(n_{2} -2) \times (n_{r}-2)}\\
     \vdots&\vdots&\vdots&\ddots&\vdots\\
     {{\mathcal{T}}_{(n_{r}-2)\times 2}}&{0}_{(n_{r}-2) \times (n_{1}-2)}&{0}_{(n_{r}-2) \times (n_{2}-2)}& \cdots&\mathcal{A}_{n_{r}-2}          
\end{array} \right],$$
where
$\mathcal{D}=[d_{ij}]_{2 \times 2}$ such that
$d_{ij}=\left \{ \begin{array}{ccccc}
   \Delta, & if & i=j&&\\
   -1,& if &{i \neq j}& and& u \sim v\\
   0, & if& i \neq j& and&u \nsim v     
\end{array} \right.$ and $\mathcal{T}=[t_{ij}]_{2 \times n_{p}-2}$ such that
$
t_{ij}=\left \{ \begin{array}{ccc}
   1, & \mbox{if } \ i=j=1, &\\
  -1, & \mbox{ if  } \ i =2 \ \mbox{ and } \ j=n_p-2,  &   \\
   0 ,& \mbox{otherwise.}
\end{array} \right. 
$ 

Observe that
$A_{n_{p}-2} \in \{ A_2, A_4, A_6\}$, where

     ${A_{2}}=\left[ \begin{array}{cc}
     2 & -1   \\
     -1 & 2    
     
\end{array} \right]$, 
${A_{4}}=\left[ \begin{array}{cccc}
     2 & -1&0&0   \\
     -1 & 2&-1&0\\
     0&-1&2&-1\\
     0&0&-1&2
     
\end{array} \right]$ and ${A_{6}}=\left[ \begin{array}{cccccc}
     2 & -1&0&0& 0&0  \\
     -1 & 2&-1&0&0&0\\
     0&1&2&1&0&0\\
     0&0&-1&2&-1&0\\
     0&0&0&-1&2&-1\\
     0&0&0&0&-1&2
     
\end{array} \right].$

By using the matrix $L(G)$ presented above, we have the following propositions.

\begin{proposition} \label{p8}
Let $G \in {\mathcal{G}^{\prime}_2}$ with $n \geq 9$ vertices. If $G$ has a subgraph $P_8$ with end vertices $u$ and $v,$ then $G$ is not $L$-integral.
\end{proposition}

\begin{proof}
Let $G \in {\mathcal{G}^{\prime}_2}$ with $n \geq 9$ vertices. Suppose that $G$ contains a path $P_8$ with the sequence of vertices $ux_1 \cdots x_6v$. Let $B_{n-2}$ be the principal submatrix of $L(G)$ obtained by removing both rows and columns that correspond to vertices $u$ and $v.$ Note that $B_{n-2}$ is a block diagonal matrix and one of its blocks is the matrix $A_6$, which has two eigenvalues in the interval $(0,1).$ From Lemma \ref{Horn}, we have $0< \theta_{n-2}(B_{n-2})<\theta_{n-3}(B_{n-2})<1$ and consequently from Proposition \ref{submatrizprincipal}, we conclude that $0<\mu_{n-1}(G)=a(G) \leq \theta_{n-3}(B_{n-2}) < 1$. Therefore, $G$ is not $L$-integral.
\end{proof}

\begin{proposition} \label{p6}
Let $G \in {\mathcal{G}^{\prime}_2}$ with $n \geq 7$ vertices. If $G$ has a subgraph $P_6$ with end vertices $u$ and $v,$ then $G$ is not $L$-integral.
\end{proposition}

\begin{proof}
Let $G \in {\mathcal{G}^{\prime}_2}$ with $n \geq 7$ vertices. Suppose that $G$ contains $a \geq 2$ paths $P_6$ with the sequence of vertices $u{x_1}^i \cdots {x_4}^iv$ for $2 \leq i \leq a$. Let $B_{n-2}$ be the principal submatrix of $L(G)$ obtained by removing both rows and columns that correspond to vertices $u$ and $v.$ Note that $B_{n-2}$ is a block diagonal matrix and $a \geq 2$ of its blocks is the matrix $A_4$. It is easy to see that $A_4$ has one eigenvalue in the interval $(0,1).$ Then, $B_{n-2}$ contains $a \geq 2$ eigenvalues in the interval $(0,1).$ From Lemma \ref{Horn}, we have $0<\theta_{n-2}(B_{n-2})<\theta_{n-3}(B_{n-2})<1,$ and consequently, from Proposition \ref{submatrizprincipal}, we conclude that $0<\mu_{n-1}(G)\leq \theta_{n-3}(B_{n-2}) < 1$. Therefore, $G$ is not $L$-integral. 

Now, suppose that $G$ contains one path $P_6$ with the sequence of vertices $ux_1 \cdots x_4v$ along with one path $P_5$ or one path $P_7$. Therefore the principal submatrix of $L(G)$ obtained by removing both rows and columns that correspond to vertices $u$ and $v$, $B_{n-2}$, is a block diagonal matrix wich has at least two blocks $A_6$ and $A_5$ or $A_6$ and $A_7$. In the both cases, $B_{n-2}$ contains $a \geq 2$ eigenvalues in the interval $(0,1)$ and, consequently, from Proposition \ref{submatrizprincipal}, we conclude that $0<\mu_{n-1}(G) \leq \theta_{n-3}(B_{n-2}) < 1$. Then, $G$ is not $L$-integral.

Finally, suppose that $G$ contains one path $P_6$, with the sequence of vertices $ux_1 \cdots x_4v$, with $s \geq 1$ paths $P_3$ and $t \geq 0$ path $P_4$. So, we need to consider the following cases:

\noindent{Case 1: $G$ contains $t\geq 1$ paths $P_4$, $s \geq 1$ paths $P_3$ and one path $P_6.$}

\noindent{Case 1.1: $u$ and $v$ are adjacent.}



By a convenient labeling for the vertices, $L(G)$ can be described the following way: 
$$L(G)=
\left[ 
\begin{array}{c|c|ccc|ccc|ccc|cccc}
    s+t+2 &-1 &-1&\cdots & -1&-1&\cdots & -1&0 &\cdots & 0 &-1 & 0 & 0 &0  \\
    \hline
     -1 &s+t+2 &-1&\cdots & -1&0&\cdots & 0&-1 &\cdots & -1& 0 & 0 & 0 &-1\\
     \hline
    -1& -1 & 2 & \cdots & 0 &0 & \cdots & 0 &0 & \cdots& 0& 0& 0& 0& 0\\
  \vdots&\vdots&\vdots&\ddots&\vdots&\vdots&\ddots&\vdots&\vdots&\ddots&\vdots&\vdots&\vdots&\vdots&\vdots\\
  -1& -1 & 0 & \cdots & 2 &0 & \cdots & 0 &0 & \cdots& 0& 0& 0& 0& 0\\
  \hline
     -1& 0 & 0 & \cdots & 0 &2 & \cdots & 0 &-1 & \cdots& 0& 0& 0& 0& 0\\
     \vdots& \vdots & \vdots & \ddots & \vdots &\vdots & \ddots & \vdots &\vdots & \ddots& \vdots& \vdots& \vdots& \vdots& \vdots \\ 
     -1&0&0&\cdots&0&0&\cdots&2&0&\cdots&-1&0&0&0&0\\ 
     \hline
      0& -1 & 0 & \cdots & 0 &-1 & \cdots & 0 &2 & \cdots& 0& 0& 0& 0& 0\\
     \vdots& \vdots & \vdots & \ddots & \vdots &\vdots & \ddots & \vdots &\vdots & \ddots& \vdots& \vdots& \vdots& \vdots& \vdots \\ 
     0&-1&0&\cdots&0&0&\cdots&-1&0&\cdots&2&0&0&0&0\\ 
     \hline
   -1&0 &0&\dots&0 &0&\dots&0&0&\dots&0&2&-1&0&0\\
   0&0 &0&\dots&0 &0&\dots&0&0&\dots&0&-1&2&-1&0\\
   0&0 &0&\dots&0 &0&\dots&0&0&\dots&0&0&-1&2&-1\\
   0&-1 &0&\dots&0 &0&\dots&0&0&\dots&0&0&0&-1&2\\
\end{array} 
\right].$$

According to Theorem \ref{equitablepartion}, the eigenvalues of the matrix

$$R_{L(G)}=\left[ \begin{array}{ccccccccc}
   s+t+2 & -1 & -s &-t & 0 & -1& 0 & 0&0 \\ 
    -1   &s+t+2& -s & 0 & -t & 0 & 0 & 0&-1 \\
    -1 & -1 & 2 & 0 & 0 & 0 & 0 & 0&0 \\ 
    -1 & 0 & 0 & 2 &  -1 & 0& 0 & 0&0 \\ 
    0 & -1 & 0 & -1 &  2 & 0& 0 & 0& 0 \\ 
    -1 & 0 & 0 & 0 &  0 & 2& -1 & 0&0  \\
    0 & 0 & 0 & 0 &  0 & -1& 2 & -1 & 0\\
    0 & 0 & 0 & 0 &  0 & 0& -1 & 2 & -1\\
    0 & -1 & 0 & 0 &  0 & 0& 0 & -1 & 2\\

\end{array} \right]$$
are eigenvalues of $L(G)$, whose characteristic polynomial is $p(\lambda)=\lambda^9+(-2s-2t-18)\lambda^8+(s^2+2st+t^2+30s+30t+137)\lambda^7+(-12s^2-24st-12t^2-184s-186t-574)\lambda^6+(56s^2+114st+57t^2+594s+614t+1443)\lambda^5+(-128s^2-272st-136t^2-1082s-1154t-2222)\lambda^4+(148s^2+338st+169t^2+1106s+1216t+2043)\lambda^3+(-80s^2-200st-100t^2-588s-650t-1026)\lambda^2+(15s^2+40st+20t^2+126s+132t+216)\lambda$. 
As $p(3)=-6st-3t^2<0$ and $p(4)=4(s+2t)(15s+10t- 18)>0$, we conclude that there is a root in the interval $(3,4)$, and consequently $G$ is not $L$-integral.

\noindent{Case 1.2: $u$ and $v$ are non-adjacent.}

By a convenient labeling for the vertices, $L(G)$ can be described the following way 
$$L(G)=
\left[ 
\begin{array}{c|c|ccc|ccc|ccc|cccc}
    s+t+1 &0 &-1&\cdots & -1&-1&\cdots & -1&0 &\cdots & 0 &-1 & 0 & 0 &0  \\
    \hline
     0 &s+t+1 &-1&\cdots & -1&0&\cdots & 0&-1 &\cdots & -1& 0 & 0 & 0 &-1\\
     \hline
    -1& -1 & 2 & \cdots & 0 &0 & \cdots & 0 &0 & \cdots& 0& 0& 0& 0& 0\\
  \vdots&\vdots&\vdots&\ddots&\vdots&\vdots&\ddots&\vdots&\vdots&\ddots&\vdots&\vdots&\vdots&\vdots&\vdots\\
  -1& -1 & 0 & \cdots & 2 &0 & \cdots & 0 &0 & \cdots& 0& 0& 0& 0& 0\\
  \hline
     -1& 0 & 0 & \cdots & 0 &2 & \cdots & 0 &-1 & \cdots& 0& 0& 0& 0& 0\\
     \vdots& \vdots & \vdots & \ddots & \vdots &\vdots & \ddots & \vdots &\vdots & \ddots& \vdots& \vdots& \vdots& \vdots& \vdots \\ 
     -1&0&0&\cdots&0&0&\cdots&2&0&\cdots&-1&0&0&0&0\\ 
     \hline
      0& -1 & 0 & \cdots & 0 &-1 & \cdots & 0 &2 & \cdots& 0& 0& 0& 0& 0\\
     \vdots& \vdots & \vdots & \ddots & \vdots &\vdots & \ddots & \vdots &\vdots & \ddots& \vdots& \vdots& \vdots& \vdots& \vdots \\ 
     0&-1&0&\cdots&0&0&\cdots&-1&0&\cdots&2&0&0&0&0\\ 
     \hline
   -1&0 &0&\dots&0 &0&\dots&0&0&\dots&0&2&-1&0&0\\
   0&0 &0&\dots&0 &0&\dots&0&0&\dots&0&-1&2&-1&0\\
   0&0 &0&\dots&0 &0&\dots&0&0&\dots&0&0&-1&2&-1\\
   0&-1 &0&\dots&0 &0&\dots&0&0&\dots&0&0&0&-1&2\\
\end{array} 
\right].$$

According to Theorem \ref{equitablepartion}, the eigenvalues of the matrix

$$R_{L(G)}=\left[ \begin{array}{ccccccccc}
   s+t+1 & 0 & -s &-t & 0 & -1& 0 & 0&0 \\ 
    0   &s+t+1& -s & 0 & -t & 0 & 0 & 0&-1 \\
    -1 & -1 & 2 & 0 & 0 & 0 & 0 & 0&0 \\ 
    -1 & 0 & 0 & 2 &  -1 & 0& 0 & 0&0 \\ 
    0 & -1 & 0 & -1 &  2 & 0& 0 & 0& 0 \\ 
    -1 & 0 & 0 & 0 &  0 & 2& -1 & 0&0  \\
    0 & 0 & 0 & 0 &  0 & -1& 2 & -1 & 0\\
    0 & 0 & 0 & 0 &  0 & 0& -1 & 2 & -1\\
    0 & -1 & 0 & 0 &  0 & 0& 0 & -1 & 2\\

\end{array} \right]$$
are eigenvalues of $L(G)$, whose characteristic polynomial is $p(\lambda)=\lambda^9+(-2s-2t-16)\lambda^8+(s^2+2st+t^2+28s+28t+107)\lambda^7+(-12s^2-24st-12t^2-160s-160t-388)\lambda^6+(56s^2+114st+57t^2+482s+480t+827)\lambda^5+(-128s^2-272st-136t^2-826s-808t-1048)\lambda^4+(148s^2+338st+169t^2+810s+754t+757)\lambda^3+(-80s^2-200st-100t^2-428s-360t-276)\lambda^2+(15^2+40st+20t^2+96s+72t+36)\lambda$. 
As $p(3)=-6st-3t^2<0$ and $p(4)=4(s+2t-2)(15s+10t-18)>0$, we conclude that there is a root in the interval $(3,4)$, and consequently $G$ is not $L$-integral.

\noindent{Case 2: $G$ contains $s \geq 1$ paths $P_3$ and one path $P_6$.}

\noindent{Case 2.1: $u$ and $v$ are adjacent.}

By a convenient labeling for the vertices, $L(G)$ can be described the following way:
$$L(G)=
\left[ 
\begin{array}{c|c|ccc|cccc}
    s+1 &-1 &-1&\cdots & -1 &-1 & 0 & 0 &0  \\
    \hline
     -1 &s+1 &-1&\cdots & -1& 0 & 0 & 0 &-1\\
     \hline
    -1& -1 & 2 & \cdots & 0 & 0& 0& 0& 0\\
  \vdots&\vdots&\vdots&\ddots&\vdots&\vdots&\vdots&\vdots&\vdots\\
  -1& -1 & 0 & \cdots & 2 & 0& 0& 0& 0\\
  \hline
   -1&0 &0&\dots&0&2&-1&0&0\\
   0&0 &0&\dots&0 &-1&2&-1&0\\
   0&0 &0&\dots&0&0&-1&2&-1\\
   0&-1 &0&\dots&0 &0&0&-1&2\\
\end{array} 
\right].$$

According to Theorem \ref{equitablepartion}, the eigenvalues of the matrix

$$R_{L(G)}=\left[ \begin{array}{ccccccc}
   s+1 & -1 & -s & -1& 0 & 0&0 \\ 
    -1   &s+1& -s &  0 & 0 & 0&-1 \\
    -1 & -1 & 2 & 0 &  0 & 0&0 \\ 
    -1 & 0 & 0 &  2& -1 & 0&0  \\
    0 & 0 & 0 &  -1& 2 & -1 & 0\\
    0 & 0 & 0 &  0& -1 & 2 & -1\\
    0 & -1 & 0 &  0& 0 & -1 & 2\\

\end{array} \right]$$
are eigenvalues of $L(G)$, whose characteristic polynomial is $p(\lambda)=\lambda^7 + (-2s - 14)\lambda^6 + (s^2 + 22s + 78)\lambda^5 + (-8s^2 - 90s - 220)\lambda^4 + (21s^2 + 168s + 329)\lambda^3 + (-20s^2 - 140s - 246)\lambda^2 + (5s^2 + 42s + 72)\lambda$. As $p(1)=-s^2<0$ and $p(2)=2s^2 + 4s>0$, we conclude that there is a root in the interval $(1,2)$, and consequently $G$ is not $L$-integral.

\noindent{Case 2.2: $u$ and $v$ are non-adjacent.}

It is easy to see that for $s=1$ or $s=2$, $G$ is not $L$-integral. Suppose that $s \geq 3$. By a convenient labeling for the vertices, $L(G)$ can be described the following way: 
$$L(G)=
\left[ 
\begin{array}{c|c|ccc|cccc}
    s+1 &0 &-1&\cdots & -1 &-1 & 0 & 0 &0  \\
    \hline
     0 &s+1 &-1&\cdots & -1& 0 & 0 & 0 &-1\\
     \hline
    -1& -1 & 2 & \cdots & 0 & 0& 0& 0& 0\\
  \vdots&\vdots&\vdots&\ddots&\vdots&\vdots&\vdots&\vdots&\vdots\\
  -1& -1 & 0 & \cdots & 2 & 0& 0& 0& 0\\
  \hline
   -1&0 &0&\dots&0&2&-1&0&0\\
   0&0 &0&\dots&0 &-1&2&-1&0\\
   0&0 &0&\dots&0&0&-1&2&-1\\
   0&-1 &0&\dots&0 &0&0&-1&2\\
\end{array} 
\right].$$

According to Theorem \ref{equitablepartion}, the eigenvalues of the matrix

$$R_{L(G)}=\left[ \begin{array}{ccccccc}
   s+1 & 0 & -s & -1& 0 & 0&0 \\ 
    0   &s+1& -s &  0 & 0 & 0&-1 \\
    -1 & -1 & 2 & 0 &  0 & 0&0 \\ 
    -1 & 0 & 0 &  2& -1 & 0&0  \\
    0 & 0 & 0 &  -1& 2 & -1 & 0\\
    0 & 0 & 0 &  0& -1 & 2 & -1\\
    0 & -1 & 0 &  0& 0 & -1 & 2\\

\end{array} \right]$$
are eigenvalues of $L(G)$, whose characteristic polynomial is $p(\lambda)=\lambda^7 + (-2s - 12)\lambda^6 + (s^2 + 20s + 56)\lambda^5 + (-8s^2 - 74s - 128)\lambda^4 + (21s^2 + 126s + 147)\lambda^3 + (-20s^2 - 100s - 76)\lambda^2 + (5s^2 + 32s + 12)\lambda$. 
As $p(2)=2s^2>0$ and $p(3)=-3s^2 + 6s<0$ for $s\geq 3$, we conclude that there is a root in the interval $(2,3)$, and consequently $G$ is not $L$-integral.
\end{proof}

\begin{proposition} \label{p4}
Let $G \in {\mathcal{G}^{\prime}_2}$ with $n \geq 6$ vertices. If $G$ has a subgraph $P_4$ with end vertices $u$ and $v,$ then $G$ is not $L$-integral.
\end{proposition}

\begin{proof}
Let $G \in {\mathcal{G}^{\prime}_2}$ with $n \geq 6$ vertices. Suppose that $G$ contains one path $P_4$ with the sequence of vertices $ux_1x_2v$ and with at least two paths in the set $\{P_5, P_7\}$. Let $B_{n-2}$ the submatrix principal of $L(G)$ obtained by removing boths rows and columns corresponding to vertices $u$ and $v$. As $B_{n-2}$ is a block diagonal matrix and its blocks belong to the set $\{A_3$, $A_5\}$ which have one eigenvalue in the interval $(0,1)$, from Proposition \ref{submatrizprincipal}, we conclude that $0<\mu_{n-1}(G)\leq \theta_{n-3}(B_{n-2}) < 1$. Then, $G$ is not $L$-integral.

Now, suppose that $G$ contains $t\geq 1$ paths $P_4$, with the sequence of vertices $ux^{q}_{1}x^{q}_{2}v$, such that $1\leq q \leq t $, $s \geq 1$ paths $P_3$ or/and one path of the set $\{P_5,P_7\}$.  So we need to analyse the following cases:

\noindent{Case 1: $G$ contains $t\geq 1$ paths $P_4$ and $s \geq 1$ paths $P_3$.}

\noindent{Case 1.1: $u$ and $v$ are adjacent.}

By a convenient labeling for the vertices, $L(G)$ can be described the following way: 
$$L(G)=
\left[ 
\begin{array}{c|c|ccc|ccc|ccc}
    s+t+1 &-1 &-1&\cdots & -1&-1&\cdots & -1&0 &\cdots & 0\\
    \hline
     -1 &s+t+1 &-1&\cdots & -1&0&\cdots & 0&-1 &\cdots & -1\\
     \hline
    -1& -1 & 2 & \cdots & 0 &0 & \cdots & 0 &0 & \cdots& 0\\
  \vdots&\vdots&\vdots&\ddots&\vdots&\vdots&\ddots&\vdots&\vdots&\ddots&\vdots\\
  -1& -1 & 0 & \cdots & 2 &0 & \cdots & 0 &0 & \cdots& 0\\
  \hline
     -1& 0 & 0 & \cdots & 0 &2 & \cdots & 0 &-1 & \cdots& 0\\
     \vdots&\vdots&\vdots&\ddots&\vdots&\vdots&\ddots&\vdots&\vdots&\ddots&\vdots\\
     -1& 0 & 0 & \cdots & 0 &0 & \cdots & 2 &0& \cdots& -1\\
     \hline
     0& -1 & 0 & \cdots & 0 &-1 & \cdots & 0 &2 & \cdots& 0\\
     \vdots&\vdots&\vdots&\ddots&\vdots&\vdots&\ddots&\vdots&\vdots&\ddots&\vdots\\
     0& -1 & 0 & \cdots & 0 &0 & \cdots & -1 &0 & \cdots& 2\\
     
\end{array} 
\right].$$

According to Theorem \ref{equitablepartion}, the eigenvalues of the matrix

$$R_{L(G)}=\left[ \begin{array}{ccccc}
   s+t+1 & -1 & -s &-t & 0 \\ 
    -1   &s+t+1& -s & 0 & -t \\
    -1 & -1 & 2 & 0 & 0  \\ 
    -1 & 0 & 0 & 2 &  -1  \\ 
    0 & -1 & 0 & -1 &  2  \\
\end{array} \right]$$
are eigenvalues of $L(G)$, whose characteristic polynomial is $p(\lambda)=\lambda^5+(-2s-2t-8)\lambda^4+(s^2+2st+t^2+12s+12t+23)\lambda^3+(-4s^2-8st-4t^2-22s-24t-28)\lambda^2+(3s^2+8st+4t^2+12s+16t+12)\lambda$. As $p(1)=2st+t^2+2t>0$ and $p(2)=-2s^2<0$, we conclude that there is a root in the interval $(1,2)$, and, consequently, $G$ is not $L$-integral.

\noindent{Case 1.2: $u$ and $v$ are non-adjacent.}

By a convenient labeling for the vertices, $L(G)$ can be described the following way: 
$$L(G)=
\left[ 
\begin{array}{c|c|ccc|ccc|ccc}
    s+t &0 &-1&\cdots & -1&-1&\cdots & -1&0 &\cdots & 0\\
    \hline
     0 &s+t &-1&\cdots & -1&0&\cdots & 0&-1 &\cdots & -1\\
     \hline
    -1& -1 & 2 & \cdots & 0 &0 & \cdots & 0 &0 & \cdots& 0\\
  \vdots&\vdots&\vdots&\ddots&\vdots&\vdots&\ddots&\vdots&\vdots&\ddots&\vdots\\
  -1& -1 & 0 & \cdots & 2 &0 & \cdots & 0 &0 & \cdots& 0\\
  \hline
     -1& 0 & 0 & \cdots & 0 &2 & \cdots & 0 &-1 & \cdots& 0\\
     \vdots&\vdots&\vdots&\ddots&\vdots&\vdots&\ddots&\vdots&\vdots&\ddots&\vdots\\
     -1& 0 & 0 & \cdots & 0 &0 & \cdots & 2 &0& \cdots& -1\\
     \hline
     0& -1 & 0 & \cdots & 0 &-1 & \cdots & 0 &2 & \cdots& 0\\
     \vdots&\vdots&\vdots&\ddots&\vdots&\vdots&\ddots&\vdots&\vdots&\ddots&\vdots\\
     0& -1 & 0 & \cdots & 0 &0 & \cdots & -1 &0 & \cdots& 2\\
     
\end{array} 
\right].$$

According to Theorem \ref{equitablepartion}, the eigenvalues of the matrix

$$R_{L(G)}=\left[ \begin{array}{ccccc}
   s+t& 0 & -s &-t & 0 \\ 
    0   &s+t& -s & 0 & -t \\
    -1 & -1 & 2 & 0 & 0  \\ 
    -1 & 0 & 0 & 2 &  -1  \\ 
    0 & -1 & 0 & -1 &  2  \\
\end{array} \right]$$
are eigenvalues of $L(G)$, whose characteristic polynomial is $p(\lambda)=\lambda^5+(-2s-2t-6)\lambda^4+(s^2+2st+t^2+10s+10t+11)\lambda^3+(-4s^2-8st-4t^2-14s-14t-6)\lambda^2+(3s^2+8st+4t^2+6s+4t)\lambda$. Then, we have:

\begin{enumerate}
    \item[{\rm (i)}] for $s=1$ and $t \geq 2$, $p(\lambda)=\lambda(\lambda^2-\lambda(4+t)+3+2t)^2$, whose roots are $0, \frac{-\sqrt{t^2 + 4} + t + 4}{2}$ with multiplicity 2, and $\frac{\sqrt{t^2 + 4} + t + 4}{2}$ with multiplicity 2 as well. As $t<\sqrt{t^2+4}<t+1$, $p(\lambda)$ has non-integer roots ;
    
    \item[{\rm (ii)}]for $s=2$ and $t=1$, $Spec_L(G)=\{4.73^{[1]}, 4^{[1]}, 2^{[2]}, 1.27^{[1]}, 0^{[1]}  \}$;
    
    \item[{\rm (iii)}]for $s=2$ and $t=2$, $Spec_L(G)=\{5.56^{[1]}, 5^{[1]}, 3^{[1]}, 2^{[2]}, 1.44^{[1]}, 1^{[1]}, 0^{[1]}  \}$;
    
    \item[{\rm (iv)}]for $s=2$ and $t \geq 3$, $p(\lambda)=(-2t - 10)\lambda^4 + \lambda^5 + (t^2 + 14t+35)\lambda^3+(-4t^2 -30t -50)\lambda^2 + (4t^2 + 20t + 24)\lambda$, whose roots are $0, 2, t+3,  \frac{-\sqrt{t^2 +2t+ 9} + t + 5}{2}, \frac{\sqrt{t^2 +2t+9} + t + 5}{2}$. As $t+1<\sqrt{t^2+2t+9}<t+2$, $p(\lambda)$ has non-integer roots ; 
    
  \item[{\rm (v)}]for $s\geq 3$ and $t \geq 1$, $p(2)=-2s^2+4s<0$ and $p(3)=6st+3t^2-6t>0$. So, we conclude that there is a root in the interval $(2,3)$.  
\end{enumerate}
Therefore, in all previous cases we obtain that $G$ is not $L$-integral.

\noindent{Case 2: $G$ contains $t\geq 1$ paths $P_4$, $s \geq 0$ paths $P_3$, and one path $P_5$.}

\noindent{Case 2.1: $u$ and $v$ are adjacent.}

By a convenient labeling for the vertices, $L(G)$ can be described in the following way: 
$$L(G)=
\left[ 
\begin{array}{c|c|ccc|ccc|ccc|ccc}
    s+t+2 &-1 &-1&\cdots & -1&-1&\cdots & -1&0 &\cdots & 0 &-1 & 0 & 0 \\
    \hline
     -1 &s+t+2 &-1&\cdots & -1&0&\cdots & 0&-1 &\cdots & -1& 0 & 0 & -1\\
     \hline
    -1& -1 & 2 & \cdots & 0 &0 & \cdots & 0 &0 & \cdots& 0& 0& 0& 0\\
  \vdots&\vdots&\vdots&\ddots&\vdots&\vdots&\ddots&\vdots&\vdots&\ddots&\vdots&\vdots&\vdots&\vdots\\
  -1& -1 & 0 & \cdots & 2 &0 & \cdots & 0 &0 & \cdots& 0& 0& 0& 0\\
  \hline
     -1& 0 & 0 & \cdots & 0 &2 & \cdots & 0 &-1 & \cdots& 0& 0& 0& 0\\
     \vdots& \vdots & \vdots & \ddots & \vdots &\vdots & \ddots & \vdots &\vdots & \ddots& \vdots& \vdots& \vdots& \vdots \\ 
     -1&0&0&\cdots&0&0&\cdots&2&0&\cdots&-1&0&0&0\\ 
     \hline
      0& -1 & 0 & \cdots & 0 &-1 & \cdots & 0 &2 & \cdots& 0& 0& 0& 0\\
     \vdots& \vdots & \vdots & \ddots & \vdots &\vdots & \ddots & \vdots &\vdots & \ddots& \vdots& \vdots& \vdots& \vdots\\ 
     0&-1&0&\cdots&0&0&\cdots&-1&0&\cdots&2&0&0&0\\ 
     \hline
   -1&0 &0&\dots&0 &0&\dots&0&0&\dots&0&2&-1&0\\
   0&0 &0&\dots&0 &0&\dots&0&0&\dots&0&-1&2&-1\\
0&-1 &0&\dots&0 &0&\dots&0&0&\dots&0&0&-1&2\\
\end{array} 
\right].$$

According to Theorem \ref{equitablepartion}, the eigenvalues of the matrix

$$R_{L(G)}=\left[ \begin{array}{cccccccc}
   s+t+2 & -1 & -s &-t & 0 & -1& 0 & 0 \\ 
    -1   &s+t+2& -s & 0 & -t & 0 &  0&-1 \\
    -1 & -1 & 2 & 0 & 0 & 0 & 0 & 0 \\ 
    -1 & 0 & 0 & 2 &  -1 & 0& 0 & 0 \\ 
    0 & -1 & 0 & -1 &  2 & 0& 0 & 0 \\ 
    -1 & 0 & 0 & 0 &  0 & 2& -1 & 0  \\
    0 & 0 & 0 & 0 &  0 & -1& 2 & -1 \\
     0 & -1 & 0 & 0 &  0 & 0 & -1 & 2\\

\end{array} \right]$$
are eigenvalues of $L(G)$, whose characteristic polynomial is $p(\lambda)=\lambda^8 + (-2s - 2t - 16)\lambda^7 + (s^2 + 2st + t^2 + 26s + 26t + 106)\lambda^6 + (-10s^2 - 20st - 10t^2 - 134s - 136t - 376)\lambda^5 + (37s^2 + 76st + 38t^2 + 348s + 364t + 770)\lambda^4 + (-62s^2 - 136st - 68t^2 - 478s - 520t - 910)\lambda^3 + (46s^2 + 112st + 56t^2 + 330s + 370t + 575)\lambda^2 + (-12s^2 - 32st - 16t^2 - 90s - 100t - 150)\lambda$. 
As $p(0.5)\approx -0.23s^2 - 0.84st - 0.42t^2 - 4.30s - 3.61t - 7.09<0$ and $p(1)=t^2 + 2t>0$, we conclude that there is a root in the interval $(0.5,1)$, and consequently $G$ is not $L$-integral.

\noindent{Case 2.2: $u$ and $v$ are non-adjacent.}

For $s=0$ and $t=2$, it is easy to see that $G$ is not $L$-integral. By a convenient labeling for the vertices, $L(G)$ can be described the following way: 
$$L(G)=
\left[ 
\begin{array}{c|c|ccc|ccc|ccc|ccc}
    s+t+1 &0 &-1&\cdots & -1&-1&\cdots & -1&0 &\cdots & 0 &-1 & 0 & 0 \\
    \hline
     0 &s+t+1 &-1&\cdots & -1&0&\cdots & 0&-1 &\cdots & -1& 0 & 0 & -1\\
     \hline
    -1& -1 & 2 & \cdots & 0 &0 & \cdots & 0 &0 & \cdots& 0& 0& 0& 0\\
  \vdots&\vdots&\vdots&\ddots&\vdots&\vdots&\ddots&\vdots&\vdots&\ddots&\vdots&\vdots&\vdots&\vdots\\
  -1& -1 & 0 & \cdots & 2 &0 & \cdots & 0 &0 & \cdots& 0& 0& 0& 0\\
  \hline
     -1& 0 & 0 & \cdots & 0 &2 & \cdots & 0 &-1 & \cdots& 0& 0& 0& 0\\
     \vdots& \vdots & \vdots & \ddots & \vdots &\vdots & \ddots & \vdots &\vdots & \ddots& \vdots& \vdots& \vdots& \vdots \\ 
     -1&0&0&\cdots&0&0&\cdots&2&0&\cdots&-1&0&0&0\\ 
     \hline
      0& -1 & 0 & \cdots & 0 &-1 & \cdots & 0 &2 & \cdots& 0& 0& 0& 0\\
     \vdots& \vdots & \vdots & \ddots & \vdots &\vdots & \ddots & \vdots &\vdots & \ddots& \vdots& \vdots& \vdots& \vdots\\ 
     0&-1&0&\cdots&0&0&\cdots&-1&0&\cdots&2&0&0&0\\ 
     \hline
   -1&0 &0&\dots&0 &0&\dots&0&0&\dots&0&2&-1&0\\
   0&0 &0&\dots&0 &0&\dots&0&0&\dots&0&-1&2&-1\\
0&-1 &0&\dots&0 &0&\dots&0&0&\dots&0&0&-1&2\\
\end{array} 
\right].$$

According to Theorem \ref{equitablepartion}, the eigenvalues of the matrix

$$R_{L(G)}=\left[ \begin{array}{cccccccc}
   s+t+1 & 0 & -s &-t & 0 & -1& 0 & 0 \\ 
    0   &s+t+1& -s & 0 & -t & 0 &  0&-1 \\
    -1 & -1 & 2 & 0 & 0 & 0 & 0 & 0 \\ 
    -1 & 0 & 0 & 2 &  -1 & 0& 0 & 0 \\ 
    0 & -1 & 0 & -1 &  2 & 0& 0 & 0 \\ 
    -1 & 0 & 0 & 0 &  0 & 2& -1 & 0  \\
    0 & 0 & 0 & 0 &  0 & -1& 2 & -1 \\
     0 & -1 & 0 & 0 &  0 & 0 & -1 & 2\\

\end{array} \right]$$
are eigenvalues of $L(G)$, whose characteristic polynomial is $p(\lambda)=\lambda^8 + (-2s - 2t - 14)\lambda^7 + (s^2 + 2st + t^2 + 24s + 24t + 80)\lambda^6 + (-10s^2 - 20st - 10t^2 - 114s - 114t - 240)\lambda^5 + (37s^2 + 76st + 38t^2 + 274s + 272t + 404)\lambda^4 + (-62s^2 - 136st - 68t^2 - 354s - 340t - 376)\lambda^3 + (46s^2 + 112st + 56t^2 + 238s + 210t + 175)\lambda^2 + (-12s^2 - 32st - 16t^2 - 66s - 52t - 30)\lambda$. 
As $p(0.5)<0$ and $p(1)>0$ for $s \neq 0$ or $t \neq 2$, we conclude that there is a root in the interval $(0.5,1)$, and consequently $G$ is not $L$-integral.

\noindent{Case 3: $G$ contains $t\geq 1$ paths $P_4$, $s \geq 1$ paths $P_3$ and one path $P_7$.}

\noindent{Case 3.1: $u$ and $v$ are adjacent.}

By a convenient labeling for the vertices, $L(G)$ can be described the following way: 
$$L(G)=
\left[ 
\begin{array}{c|c|ccc|ccc|ccc|ccccc}
    s+t+2 &-1 &-1&\cdots & -1&-1&\cdots & -1&0 &\cdots & 0 &-1 & 0 & 0& 0& 0 \\
    \hline
     -1 &s+t+2 &-1&\cdots & -1&0&\cdots & 0&-1 &\cdots & -1& 0 & 0 &0&0&  -1\\
     \hline
    -1& -1 & 2 & \cdots & 0 &0 & \cdots & 0 &0 & \cdots& 0& 0& 0& 0& 0& 0\\
  \vdots&\vdots&\vdots&\ddots&\vdots&\vdots&\ddots&\vdots&\vdots&\ddots&\vdots&\vdots&\vdots&\vdots&\vdots&\vdots\\
  -1& -1 & 0 & \cdots & 2 &0 & \cdots & 0 &0 & \cdots& 0& 0& 0& 0& 0& 0\\
  \hline
     -1& 0 & 0 & \cdots & 0 &2 & \cdots & 0 &-1 & \cdots& 0& 0& 0& 0& 0& 0\\
     \vdots& \vdots & \vdots & \ddots & \vdots &\vdots & \ddots & \vdots &\vdots & \ddots& \vdots& \vdots& \vdots& \vdots& \vdots& \vdots \\ 
     -1&0&0&\cdots&0&0&\cdots&2&0&\cdots&-1&0&0&0&0&0\\ 
     \hline
      0& -1 & 0 & \cdots & 0 &-1 & \cdots & 0 &2 & \cdots& 0& 0& 0& 0&0&0\\
     \vdots& \vdots & \vdots & \ddots & \vdots &\vdots & \ddots & \vdots &\vdots & \ddots& \vdots& \vdots& \vdots& \vdots& \vdots& \vdots\\ 
     0&-1&0&\cdots&0&0&\cdots&-1&0&\cdots&2&0&0&0&0&0\\ 
     \hline
   -1&0 &0&\dots&0 &0&\dots&0&0&\dots&0&2&-1&0&0&0\\
   0&0 &0&\dots&0 &0&\dots&0&0&\dots&0&-1&2&-1&0&0\\
   0&0 &0&\dots&0 &0&\dots&0&0&\dots&0&0&-1&2&-1&0\\
   0&0 &0&\dots&0 &0&\dots&0&0&\dots&0&0&0&-1&2&-1\\
0&-1 &0&\dots&0 &0&\dots&0&0&\dots&0&0&0&0&-1&2\\
\end{array} 
\right].$$

According to Theorem \ref{equitablepartion}, the eigenvalues of the matrix

$$R_{L(G)}=\left[ \begin{array}{cccccccccc}
   s+t+2 & -1 & -s &-t & 0 & -1& 0&  0 & 0& 0 \\ 
    -1   &s+t+2& -s & 0 & -t & 0 &0&0&     0&-1 \\
    -1 & -1 & 2 & 0 & 0 & 0 & 0 & 0&0& 0 \\ 
    -1 & 0 & 0 & 2 &  -1 & 0& 0 & 0&0& 0 \\ 
    0 & -1 & 0 & -1 &  2 & 0& 0 & 0&0& 0 \\ 
    -1 & 0 & 0 & 0 &  0 & 2& -1 & 0&0 & 0 \\
    0 & 0 & 0 & 0 &  0 & -1& 2 & -1&0 & 0\\
    0 & 0 & 0 & 0 &  0 &0 & -1& 2 & -1 & 0\\
    0 & 0 & 0 & 0 &  0 &0 &0 & -1& 2 & -1 \\
     0 & -1 & 0 & 0 &  0 &0 & 0 &0 & -1 & 2\\

\end{array} \right]$$
are eigenvalues of $L(G)$, whose characteristic polynomial is $p(\lambda)=\lambda^{10} + (-2s - 2t - 20)\lambda^9 + (s^2 + 2st + t^2 + 34s + 34t + 172)\lambda^8 + (-14s^2 - 28st - 14t^2 - 242s - 244t - 832)x^7 + (79s^2 + 160st + 80t^2 + 936s + 960t + 2485)\lambda^6 + (-230s^2 - 480st - 240t^2 - 2136s - 2246t - 4732)\lambda^5 + (367s^2 + 806st + 403t^2 + 2922s + 3160t + 5719)\lambda^4 + (-314s^2 - 740st - 370t^2 - 2324s - 2562t - 4214)\lambda^3 + (129s^2 + 328st + 164t^2 + 980s + 1064t + 1715)\lambda^2 + (-18s^2 - 48st - 24t^2 - 168s - 168t - 294)\lambda$. 
As $p(0.5)\approx 0.88s^2 + 3.16st + 1.58t^2 - 0.76s + 2.29t - 2.47>0$ and $p(1)=-4t<0$, we conclude that there is a root in the interval $(0.5,1)$, and consequently $G$ is not $L$-integral.

\noindent{Case 3.2: $u$ and $v$ are not adjacent.}

By a convenient labeling for the vertices, $L(G)$ can be described the following way: 
$$L(G)=
\left[ 
\begin{array}{c|c|ccc|ccc|ccc|ccccc}
    s+t+1 &0 &-1&\cdots & -1&-1&\cdots & -1&0 &\cdots & 0 &-1 & 0 & 0& 0& 0 \\
    \hline
     0 &s+t+1 &-1&\cdots & -1&0&\cdots & 0&-1 &\cdots & -1& 0 & 0 &0&0&  -1\\
     \hline
    -1& -1 & 2 & \cdots & 0 &0 & \cdots & 0 &0 & \cdots& 0& 0& 0& 0& 0& 0\\
  \vdots&\vdots&\vdots&\ddots&\vdots&\vdots&\ddots&\vdots&\vdots&\ddots&\vdots&\vdots&\vdots&\vdots&\vdots&\vdots\\
  -1& -1 & 0 & \cdots & 2 &0 & \cdots & 0 &0 & \cdots& 0& 0& 0& 0& 0& 0\\
  \hline
     -1& 0 & 0 & \cdots & 0 &2 & \cdots & 0 &-1 & \cdots& 0& 0& 0& 0& 0& 0\\
     \vdots& \vdots & \vdots & \ddots & \vdots &\vdots & \ddots & \vdots &\vdots & \ddots& \vdots& \vdots& \vdots& \vdots& \vdots& \vdots \\ 
     -1&0&0&\cdots&0&0&\cdots&2&0&\cdots&-1&0&0&0&0&0\\ 
     \hline
      0& -1 & 0 & \cdots & 0 &-1 & \cdots & 0 &2 & \cdots& 0& 0& 0& 0&0&0\\
     \vdots& \vdots & \vdots & \ddots & \vdots &\vdots & \ddots & \vdots &\vdots & \ddots& \vdots& \vdots& \vdots& \vdots& \vdots& \vdots\\ 
     0&-1&0&\cdots&0&0&\cdots&-1&0&\cdots&2&0&0&0&0&0\\ 
     \hline
   -1&0 &0&\dots&0 &0&\dots&0&0&\dots&0&2&-1&0&0&0\\
   0&0 &0&\dots&0 &0&\dots&0&0&\dots&0&-1&2&-1&0&0\\
   0&0 &0&\dots&0 &0&\dots&0&0&\dots&0&0&-1&2&-1&0\\
   0&0 &0&\dots&0 &0&\dots&0&0&\dots&0&0&0&-1&2&-1\\
0&-1 &0&\dots&0 &0&\dots&0&0&\dots&0&0&0&0&-1&2\\
\end{array} 
\right].$$

According to Theorem \ref{equitablepartion}, the eigenvalues of the matrix

$$R_{L(G)}=\left[ \begin{array}{cccccccccc}
   s+t+1 & 0 & -s &-t & 0 & -1& 0&  0 & 0& 0 \\ 
    0   &s+t+1& -s & 0 & -t & 0 &0&0&     0&-1 \\
    -1 & -1 & 2 & 0 & 0 & 0 & 0 & 0&0& 0 \\ 
    -1 & 0 & 0 & 2 &  -1 & 0& 0 & 0&0& 0 \\ 
    0 & -1 & 0 & -1 &  2 & 0& 0 & 0&0& 0 \\ 
    -1 & 0 & 0 & 0 &  0 & 2& -1 & 0&0 & 0 \\
    0 & 0 & 0 & 0 &  0 & -1& 2 & -1&0 & 0\\
    0 & 0 & 0 & 0 &  0 &0 & -1& 2 & -1 & 0\\
    0 & 0 & 0 & 0 &  0 &0 &0 & -1& 2 & -1 \\
     0 & -1 & 0 & 0 &  0 &0 & 0 &0 & -1 & 2\\

\end{array} \right]$$
are eigenvalues of $L(G)$, whose characteristic polynomial is $p(\lambda)=\lambda^{10} + (-2s - 2t - 18)\lambda^9 + (s^2 + 2st + t^2 + 32s + 32t + 138)\lambda^8 + (-14s^2 - 28st - 14t^2 - 214s - 214t - 588)\lambda^7 + (79s^2 + 160st + 80t^2 + 778s + 776t + 1523)\lambda^6 + (-230s^2 - 480st - 240t^2 - 1676s - 1654t - 2462)\lambda^5 + (367s^2 + 806st + 403t^2 + 2188s + 2098t + 2449)\lambda^4 + (-314s^2 - 740st - 370t^2 - 1696s - 1528t - 1414)\lambda^3 + (129s^2 + 328st + 164t^2 + 722s + 584t + 413)\lambda^2 + (-18s^2 - 48st - 24t^2 - 132s - 96t - 42)\lambda$. 
As $p(0.5)\approx 0.88s^2 + 3.16st + 1.58t^2 - 2.52s - 2.99t + 1.33>0$ and $p(1)=-4t<0$, we conclude that there is a root in the interval $(0.5,1)$, and consequently $G$ is not $L$-integral.
\end{proof}

\begin{thm} \label{teoremagrafonaobipartidofinal}
Let  $G \in {\mathcal{G}^{\prime}_2}$ with $n \geq 9$ vertices. Then $G$ is $L$-integral if  and only if $G \cong K_1 \vee (r\cdot K_1 \cup s\cdot K_2 \cup K_{1,t})$,  where $t \geq 2$ and $r+s \geq 2$ or $G \cong K_2 \vee (n-2) \cdot K_1$. 
\end{thm}

\begin{proof}
Let  $G \in {\mathcal{G}^{\prime}_2}$ with $n \geq 9$ vertices. Suppose that $G$  is $L$-integral. From Theorem \ref{a=k1} and Propositions \ref{maiorque p9}, \ref{p5+p7}, \ref{P3+P5}, \ref{P3+P7}, \ref{p8}, \ref{p6} and \ref{p4}  we conclude that $G \cong K_1 \vee (r\cdot K_1 \cup s\cdot K_2 \cup K_{1,t})$,  where $t \geq 2$ and $r+s \geq 2$ or $G \cong K_2 \vee (n-2)\cdot K_1$ and the result follows.
\end{proof}

From Theorems \ref{qintegral}, \ref{teoremaG1}  and \ref{teoremagrafonaobipartidofinal}, the proof of Theorem 1.1 is complete.



{}


\begin{thebibliography}{}






\bibitem{bro} A.E. Brouwer,  W.H. Haemers,  Spectra of Graphs,   Springer, New York, 2012.



\bibitem{fie}  M. Fiedler, Algebraic Connectivity of Graphs, Czechoslovak Mathematical Journal, 23 (1973) 298-305.

\bibitem{heu} J. van den Heuvel, Hamilton cycles and eigenvalues of graphs, Linear Algebra and Its Applications, 226-228:723-730, (1995).

\bibitem{horn} R.A. Horn, C.R. Johnson, Matrix Analysis, Cambridge University Press, 2 edition, (2013).

\bibitem{hof} P. Hof, D. Paulusma, A new characterization of $P_6$-free graphs, Discrete Applied Mathematics, 158 (2010) 731-740.

\bibitem{kir} S.J. Kirkland, J.J. Molitierno, M. Neumann, B.L. Shader, On graphs with equal algebraic and vertex connectivity, Linear Algebra and its Applications, 341 (2002) 45-56.

\bibitem{kir2} S.J. Kirkland, Constructably Laplacian integral graphs, Linear Algebra and its Applications, 423 (2007) 3-21.

\bibitem{kir3} S.J. Kirkland, Laplacian Integral Graphs with Maximum Degree 3, The Electronic Journal of Combinatorics 15 (2008) R120.



\bibitem{oli1} L.S. de Lima, N.M. de Abreu, C.S. Oliveira, Laplacian Integral graphs in S(a,b), Linear Algebra and its Applications, 423 (2007) 136-145.


\bibitem{mer1} R. Merris, Degree Maximal graphs are Laplacian integral, Linear Algebra and Its Applications, 199 (1994) 381-389.

\bibitem{mer} R. Merris, Laplacian Graph Eigenvectores, Linear Algebra and Its Applications, 278 (1998) 221-236.  

\bibitem{nov} A.F. Novanta, L.S. de Lima, C.S. Oliveira, $Q$-integral graphs with at most two vertices of degree greater than or equal to three. Available from  \mbox{https://doi.org/10.1016/j.laa.2020.03.027.}  






\end{thebibliography}
\end{document}